\documentclass[a4paper,draft,leqno,12pt]{article}

\usepackage{amsmath}
\usepackage{amscd}
\usepackage{amssymb}
\usepackage{enumerate}
\usepackage{indentfirst}
\usepackage{latexsym}
\usepackage{multicol}
\usepackage{pst-fill,pstricks}
\usepackage{theorem}
\numberwithin{equation}{section}
\def\div{\, \hbox{\rm div}\,  }

\newcommand{\R}{{\mathbb R}}
\newcommand{\N}{{\mathbb N}}

\newcommand{\Z}{{\mathbb Z}}
\newcommand{\T}{{\mathbb T}}

\def\d{\partial}

\begin{document}

	\centerline{\large Transport equation: extension of classical results for div $b\in$ BMO}

\bigskip

\begin{center}

\textbf{Piotr Bogus\l aw Mucha } \\
\text{ \ }\\
\text{Instytut Matematyki Stosowanej i Mechaniki} \\
\text{Uniwersystet Warszawski}\\
\text{ul. Banacha 2, Warszawa, Poland}\\
\text{email:} \texttt{p.mucha@mimuw.edu.pl}

\end{center}

\bigskip

{\bf Abstract.} We investigate the transport equation: $u_t+b \cdot \nabla u=0$. Our result improves
the criteria on uniqueness  of weak solutions, replacing the classical condition:
  $\div b \in L_\infty$ by $\div b \in BMO$.

\medskip



\section{Introduction}

The goal of this note is to improve  classical results concerning the Cauchy problem for the transport equation. The basis of our analysis is the following system
\begin{equation}\label{i1}
 \begin{array}{lcr}
 ×u_t+b\cdot \nabla u =0 & {\rm in} & \R^n \times (0,T),\\
u|_{t=0}=u_0 & {\rm on} & \R^n,
\end{array}
\end{equation}
where $u$ is an unknown scalar function, $b$ -- some given vector field and $u_0$ is an initial datum.

The transport equation  is one of the most fundamental examples  in the theory of partial differential equations. It describes the motion of matter under influence of the velocity
field $b$. Classically, for smooth data $b$ and $u_0$, (\ref{i1}) is solvable elementary
 by the method of characteristics. In the language of the fluid mechanics,
(\ref{i1}) says that $u$ is constant along streamlines defined by 
 the Lagrangian coordinates. This physical interpretation gives enough 
 reasons for  (\ref{i1}) to be
intensively studied from the mathematical point of view. Here we want to concentrate on 
 the optimal/critical regularity of the vector  field $b$ to control 
the existence and uniqueness of the solutions.

The classical results require that the vector field  $b$ must fulfill
\begin{equation}\label{i2}
 \div b \in  L_1(0,T,L_\infty(\R^n)),
\end{equation}
then we are able to obtain existence and uniqueness of weak solutions to (\ref{i1}). The objective
of our investigations is to relax the condition (\ref{i2}) to
\begin{equation}\label{i3}
 \div b \in L_1(0,T;BMO(\R^n)).
\end{equation}
Let us observe that this ``slightly'' broader class than  (\ref{i2}) is on the boundary of known counterexamples \cite{DL}. For any $p < \infty$ we are able to construct such $b \in W^1_p(\R^n)$ 
(time independent) to obtain an example of
the loss of uniqueness to (\ref{i1}). On the other hand the $BMO$-space appears naturally
in many considerations, since it is the limit space for the embedding $W^1_n(\R^n) \subset
BMO(\R^n)$, where  the $L_\infty$-space is not reached. 
We are able to  prove existence and uniqueness of weak solutions to (\ref{i1}) in the case of bounded solutions and improve the
uniqueness criteria for $L_p$-solutions.  Additionally we show a result concerning stability with respect to initial data.  Our approach  follows from techniques introduced in \cite{MR} to improve
the uniqueness criteria for the Euler system in bounded domains.

Results  fundamental  for our issue have been stated by R.J. DiPerna and P.L. Lions in \cite{DL}, where general questions concerning the well posedness of the problem found positive answer under condition (\ref{i2}). An interesting extension 
of the theory has been made by L. Amrozio \cite{A}, for the case of bounded solutions
replacing the condition $b \in W^1_1(\R^n)$ by $b\in BV(\R^n)$.
In the literature one can find also numerous works on generalizations of the mentioned results on broader class of function spaces \cite{AC},\cite{D2},\cite{HLL},\cite{LL1},\cite{LL2}, but positive answers  still require   condition (\ref{i2}).

In the present note we consider weak solutions meant in the following sense:

\smallskip

We say that $u \in L_\infty(0,T;L_p(\R^n))$  is a weak solution to \eqref{i1} iff the following integral identity holds 
\begin{equation}\label{w2int}
\int_0^T \int_{\R^n} u \phi_t dx dt + \int_0^T \int_{\R^n} \div\, b \, u \phi dx
+ \int_0^T \int_{\R^n} b \cdot \nabla \phi u dx dt =
-\int_{\R^n}  u_0 \phi (\cdot,0)dx
\end{equation}
for each $\phi \in C^{\infty}([0,T];C_0^\infty(\R^n))$ such that $\phi|_{t=T}\equiv 0$.

\smallskip

Let us state the main results of this paper.  First we start with the case of pointwise bounded solutions, in that case our technique
delivers  the most complete result.

\smallskip

{\it {\bf Theorem A.} Let $T>0$, $b\in L_1(0,T;W^1_{1(loc)}(\R^n))$, $u_0 \in L_\infty (\R^n)$, additionally we assume
\begin{equation}\label{a1}
\div\, b \in L_1(0,T;BMO(\R^n)), \qquad \frac{b}{1+|x|} \in L_1(0,T;L_1(\R^n))
\end{equation}
\begin{equation}\label{a2}
\mbox{ ~~~ and ~~ supp } \div\,b(\cdot,t) \subset B(0,R) \mbox{~~for a fixed $R>0$,}
\end{equation}
where $B(0,R)$ denotes the ball  centered at the origin with radius $R$.

 Then there exists a unique weak solution to the system (\ref{i1})
such that 
\begin{equation}\label{a3}
u\in L_\infty(0,T;L_\infty(\R^n)) \mbox{~~~ and ~~~}
 \|u\|_{L_\infty(0,T;L_\infty(\R^n))}\leq \|u_0\|_{L_\infty(\R^n)}.
\end{equation}}

The above result guarantees  not only the uniqueness of solutions, but also their existence. It is a consequence of a maximum principle, which is valid for the $L_\infty$-solutions. To show (\ref{a3})  condition (\ref{a1}) is not needed. The main difference to the  classical results \cite{DL} is that having (\ref{i2}) we are able to construct the $L_p$-estimates of the solutions
for finite $p$. In our case  the condition (\ref{i3}) is too weak to get such information. Additionally we are required to add an extra condition  \eqref{a2}, which is the price of our improvement of this classical criteria.

The next result concerns  stability of  solutions obtained in  Theorem A with respect to perturbations of  initial data in lower spaces.

\smallskip

{\it {\bf Theorem B.} Let $1 \leq p <\infty$ and  $b$ fulfill assumptions of Theorem A. Let $u_0,u_0^k \in L_\infty(\R^n)$
and $(u_0^k-u_0) \in L_p(\R^n)$ such that
	$\sup_{k \in \N} \|u_0^k\|_{L_\infty(\R^n)}+\|u_0\|_{L_\infty(\R^n)} =m < \infty$
and
$(u_0^k-u_0) \to 0 \mbox{~~in~} L_p(\R^n) \mbox{~~as~~} k \to \infty.$
Then 
\begin{equation}\label{b3}
(u^k-u) \to 0 \mbox{~~in~} L_\infty(0,T;L_p(\R^n)) \mbox{~~as~~} k \to \infty.
\end{equation}}

\smallskip

The last result concerns the uniqueness criteria for  $L_p$-solutions to (\ref{i1}).

\smallskip

{\it {\bf Theorem C.} Let $1 < p < \infty$, $\frac 1p + \frac{1}{p'}=1$, $b \in L_1(0,T;W^1_{p'(loc)}(\R^n))$ and conditions (\ref{a1}), (\ref{a2}) be fulfilled.
Let $u^1,u^2$ be two weak solutions to (\ref{i1}) with the same initial datum and
$u^1,u^2 \in L_\infty(0,T;L_p(\R^n))$; then $u^1 \equiv u^2$.
}

\smallskip

The above three results are proved by a reduction of considerations to an ordinary differential inequality of 
the form
\begin{equation}\label{cc4}
 \dot x=x \ln x, \qquad x|_{t=0}=0.
\end{equation}
The Osgood lemma yields the uniqueness to the system \eqref{cc4}.
This observation forms our chain of estimations to have a possibility to adapt information obtained from the Gronwall inequality. Due to low regularity of solutions, our analysis requires a special approach.
The main tool, which enables us to show the main inequality in the form of \eqref{cc4},
   is  Theorem D stated  below. 

\smallskip

{\it {\bf Theorem D.} Let $f \in BMO(\R^n)$, the support of  $f$ be bounded in $\R^n$ and
$g \in L_1(\R^n)\cap L_\infty(\R^n)$, then 
\begin{equation}\label{d1}
| \int_{\R^n} fg dx| \leq C_0\|f\|_{BMO(\R^n)}\|g\|_{L_1(\R^n)}\left[ | \ln \|g\|_{L_1(\R^n)}| + \ln (e+\|g\|_{L_\infty(\R^n)})
\right].
\end{equation}}

\smallskip

The above inequality  can be viewed as a representant of the family of logarithmic Sobolev inequalities \cite{BKM},\cite{Da1},\cite{KOZ}, however there is one important difference between  this one and others, here an extra information about derivatives of the function is not required, 
in contrast to  $L_\infty-BMO$ inequalities.
The crucial assumption is the boundedness of the
support of the function $f$, it is a consequence of results of the classical theory \cite{Tor},\cite{Z} and unfortunately it is not expected that it could be possible to omit this restriction. Methods of proving 
 \eqref{d1}  distinguish this result from others, too. They base on relations between the Zygmund space $L \ln L$ and Riesz operators.
Theorem D has been  proved in \cite{MR}, applied successfully to the evolutionary Euler system.
Outlines of the proof of Theorem D one can find in the Appendix.

The below remark shows us  a possible generalization of  stated  theorems.

\smallskip

{\it
{\bf Remark.} The results stated in Theorems A, B and C can be easily extended on  the following linear system
\begin{equation}\label{ee1}
 \begin{array}{lcr}
 ×
u_t+b\cdot \nabla u =c u + f & \mbox{in} & \Omega \times (0,T),\\
u|_{t=0}=u_0 & \mbox{on} & \Omega
\end{array}
\end{equation}
in an arbitrary domain $\Omega \subset \R^n$ with a sufficiently smooth 
boundary $\partial \Omega$, say $C^{0+1}$, enough to allow  integration by parts, and with given
$$
c, f \in L_1(0,T;L_\infty(\R^n)) \mbox{~~~and ~~~} b\cdot n =0 \mbox{ ~~ on ~~ } \partial \Omega \times (0,T),
$$
where $n$ is the normal vector to the boundary $\partial \Omega$.

Additionally, we find a natural generalization of (\ref{a1})-(\ref{a2})
\begin{equation}\label{ee2}
\begin{array}{c}
\displaystyle
 \div\,b= H_\infty+ \sum_{k=1}^\infty H_k \mbox{ such that } \\
\displaystyle
H_\infty \in L_1(0,T;L_\infty(\Omega)), \qquad H_k \in L_1(0,T;BMO(\Omega)) \mbox{ and } \\
\displaystyle
\sum_{k=1}^\infty \|H_k\|_{L_1(0,T;BMO(\Omega))} < \infty \mbox{ ~~~ with ~~~ } \sup_{k�\in \N} diam \; supp\; H_k < \infty.
\end{array}
\end{equation}
}

In the case of  bounded $\Omega$ condition \eqref{ee2} simplifies 
itself and \eqref{a2} is automatically satisfied.
We leave  the proof of Remark to a kind reader, it is almost the same as for (\ref{i1}), 
the estimations are just more technical, but the core of the problem is the same.

\smallskip

Thoughout the paper we use  standard notation. $L_p(\R^n)$ denotes the common Lebesgue space, generic constant are denoted by $C$.  Let us recall only the definition of the BMO-space.
We say that  $f \in BMO(\R^n)$, if  $f$ is locally integrable and the corresponding  semi-norm
\begin{equation}\label{bmo}
\|f\|_{BMO(\R^n)} = \sup_{x \in \R^n, r>0} 
\int_{B(x,r)} |f(y)-\{f\}_{B(x,r)}| dy
\end{equation}
is finite, where
$\{f\}_{B(x,r)}=\frac{1}{|B(x,r)|}   \int_{B(x,r)} f(y)dy$ and $B(x,r)$ is a ball with radius $r$  centered at $x$ -- see \cite{S}.
The above definition implies that (\ref{bmo}) is a semi-norm only, however in our case from assumptions on
$\div b$ follows
$$
\| \div b \|_{L_1(\R^n)} \leq | supp \div b | \| \div b \|_{BMO(\R^n)}
$$
which is a consequence of the properties of the support restricted by (\ref{a2}).

\section{Proof of Theorem A}

Our first aim is to prove the existence of weak solutions to (\ref{i1}).
To construct them we find the following sequence of
approximation of the function $b$ and initial datum $u_0$.
We require that
$b^\epsilon \in C^\infty(\R^n\times (0,T)) \mbox{~~~and~~~}
supp\, \div\,b^\epsilon(\cdot,t) \subset B(0,2R)$
and $b^\epsilon \to b \mbox{~~in~~} L_1(0,T;W^1_{1(loc)}(\R^n))$ with
suitable behavior of the norms.
For the initial datum we find
$u_0^\epsilon \in C^\infty_0(\R^n)$ with
$u_0^\epsilon \rightharpoonup^* u_0$ in $L_\infty(\R^n)$ as $\epsilon \to 0$ and
\begin{equation}\label{a6}
\|u_0^\epsilon\|_{L_\infty(\R^n)}\leq \|u_0\|_{L_\infty(\R^n)}.
\end{equation}
Then we consider the following equation with smooth coefficients $b^\epsilon$ and initial data $u_0^\epsilon$.
\begin{equation}\label{eqe}
 \begin{array}{lcr}
 �u^\epsilon_t+b^\epsilon\cdot \nabla u^\epsilon =0 & {\rm  in} & \R^n \times (0,T),\\
u^\epsilon|_{t=0}=u^\epsilon_0 & {\rm on} & \R^n.
\end{array}
\end{equation}

The method of characteristic implies the existence of smooth solutions to (\ref{eqe}) for $t\in (0,T)$. Omitting the characteristic coordinates system
we find an $\epsilon$-independent estimate for the solutions by a simple Mouser's technique. For any even $p$ we easily get
\begin{equation*}
\frac 1p \frac{d}{dt} \int {u^\epsilon}^p dx -\frac{1}{p} \int \div\, b^\epsilon \, {u^\epsilon}^p dx=0,
\end{equation*}
which implies
\begin{equation*}
\|u^\epsilon\|_{L_p(\R^n)}\leq \|u_0^\epsilon\|_{L_p(\R^n)} e^{\frac 1p \int_0^t \|\div\,b^\epsilon\|_{L_\infty(\R^n)} ds}.
\end{equation*}
Since $\div \,u^\epsilon$ is smooth and $\epsilon>0$ is fixed, hence letting $p\to \infty$, we get immediately
\begin{equation}\label{e1}
\| u^\epsilon \|_{L_\infty(0,T;L_\infty(\R^n))}\leq \|u_0^\epsilon\|_{L_\infty(\R^n)}.
\end{equation}
Note that we do not use any uniform bound on $\div b^\epsilon$.

Now we pass to the limit with $\epsilon \to 0$ in \eqref{eqe}. The  solutions to (\ref{eqe}) are classical, in particular it implies they fulfill the following integral identity
\begin{equation}\label{w1}
-\int_0^T \int_{\R^n} u^\epsilon \phi_t dx dt - \int_0^T \int_{\R^n} \div\, b^\epsilon u^\epsilon \phi dx
- \int_0^T \int_{\R^n} b^\epsilon \cdot \nabla \phi u^\epsilon dx dt =
\int_{\R^n} u_0^\epsilon \phi(\cdot,0)dx
\end{equation}
for any $\phi \in C^\infty([0,T]; C^\infty_0(\R^n))$ such that $\phi|_{t=T}\equiv 0$.

Estimates \eqref{a6} and  (\ref{e1}) imply that for a subsequence $\epsilon_k\to 0$
\begin{equation}\label{w1a}
u^{\epsilon_k} \rightharpoonup^* u \mbox{~~in~} L_\infty( 0,T;L_\infty(\R^n)) 
\mbox{~~with~~} \|u\|_{L_\infty(0,T;L_\infty(\R^n))} \leq \|u_0\|_{L_\infty(\R^n))}.
\end{equation}
Then taking the limit of (\ref{eqe}) for $\epsilon_k\to 0$, by properties of sequences
$\{b^\epsilon\}$ and $\{u_0^\epsilon\}$,  we obtain
\begin{equation}\label{w2}
-\int_0^T \int_{\R^n} u \phi_t dx dt - \int_0^T \int_{\R^n} \div\, b \, u \phi dx
- \int_0^T \int_{\R^n} b \cdot \nabla \phi u dx dt =
\int_{\R^n}  u_0 \phi (\cdot,0)dx
\end{equation}
for the same set of  test functions as in (\ref{w1}). To simplify the notation we will write $\epsilon \to 0$ instead of $\epsilon_k \to 0$.

This way the limit $u$, defined by (\ref{w1a}), is a weak solution to (\ref{i1}). It is clear that (\ref{w2}) it is its distributional version
-- \eqref{w2int}.
However the high regularity of  test functions required in (\ref{w2}) does not allow us to obtain any information concerning the uniqueness of solutions to \eqref{w2} in a direct way. 
To solve this issue we start with an application of the standard procedure. We introduce 
\begin{equation}\label{w3a}
S_\epsilon(f)=m_\epsilon \ast f= \int_{\R^n} m_\epsilon(\cdot -y)f(y) dy,
\end{equation}
where $m_\epsilon$ is a smooth function with suitable properties tending weakly to the Dirac delta -- see  (\ref{x1c}) in the Appendix.
Applying the above operator to (\ref{i1}) we get
\begin{equation}\label{w3b}
\partial_t S_\epsilon(u)+S_\epsilon(b\cdot \nabla u)=0,
\end{equation}
where $b\cdot \nabla u=\div(bu)- u\,\div\,b$ and  the r.h.s. is well defined as a distribution. In fact (\ref{w3b}) implies that
$\partial_t S_\epsilon(u)$ is well defined as a Lebesgue function, too.

We state  equation \eqref{w3b} as follows
\begin{equation}\label{w4}
\partial_t S_\epsilon(u)+b\cdot \nabla S_\epsilon(u)=R_\epsilon,
\mbox{ ~~ where ~~ } 
R_\epsilon=b\cdot \nabla S_\epsilon(u) - S_\epsilon(b\cdot \nabla u).
\end{equation}
Standard facts follow (see \eqref{x2} in Appendix) the remeider is controlled in the limit:
$R_\epsilon \to 0 \mbox{~ in ~} L_1(0,T;L_{1(loc)}(\R^n)).$
Since $R_\epsilon$ convergences locally in space, only, we introduce  a smooth function $\pi_r: \R^n \to
[0,1]$ such that $\pi_r(x)=\pi_1(\frac xr) \mbox{~~and~~}$
\begin{equation}\label{w6}
\pi_1(x)=\left\{
\begin{array}{cc}
 �1 & |x| < 1 \\
\in [0,1] & 1\leq |x| \leq 2 \\
0 & |x| > 2
\end{array}
\right.
\mbox{ ~~~~~~with ~~ } |\nabla \pi_r|\leq \frac Cr.
\end{equation}

In order to prove the uniqueness for  our system it is enough to consider (\ref{w4}) with zero initial data (due to its linearity).
Since we are forced to localize the problem, we multiply \eqref{w4} by
$S_\epsilon(u) \pi_r$ and integrate over the space, getting
\begin{equation}\label{w7}
\begin{array}{c}
\displaystyle
\frac 12 \frac{d}{dt} \int_{\R^n} (S_\epsilon(u))^2\pi_r dx - \frac 12 \int_{\R^n} \div \,b \,  (S_\epsilon(u))^2 \pi_rdx
-\frac 12 \int_{\R^n} b \cdot \nabla \pi_r \, (S_\epsilon(u))^2 dx\\
\displaystyle
=\int_{\R^n} R_\epsilon \, S_\epsilon(u)\pi_rdx.
\end{array}
\end{equation}
Then integrating over time, using properties of $S_\epsilon$ and letting $\epsilon \to 0$, next   
 differentiating with respect  $t$ we obtain
\begin{equation}\label{w8}
\frac 12 \frac{d}{dt} \int_{\R^n} u^2 \pi_r dx - \frac 12 \int_{\R^n} \div \,b\,u^2 \pi_r dx=\frac 12 \int_{\R^n} b \cdot \nabla \pi_r u^2 dx
\mbox{~~for~~} r>0.
\end{equation}
The r.h.s. of  (\ref{w8})  is estimated as follows
\begin{equation}\label{w9}
| \int_{\R^n} b \cdot \nabla \pi_r \, u^2 dx| \leq C \|u\|^2_{L_\infty(\R^n)}
\int_{\R^n \setminus B(0,r)} \frac{|b|}{1+|x|} (1+|x|)|\nabla \pi_r| dx \to 0
\mbox{~~~as~~} r \to \infty.
\end{equation}
By  definition  $(1+|x|)|\nabla \pi_r| \leq C$, because  the support of $\nabla \pi_k$
is a subset of the set: $\{ r\leq |x| \leq 2r\}$. By (\ref{w1a}) the norm $\|u\|_{L_\infty}$
is controlled, too.
Then letting  $r \to \infty$ in (\ref{w9}), we get
\begin{equation}\label{w10}
\frac 12 \frac{d}{dt} \int_{\R^n} u^2  dx - \frac 12 \int_{\R^n} \div \,b \, u^2  dx=0
\mbox{~~~~~with ~~ $u|_{t=0}\equiv 0$.}
\end{equation}

On the other hand, the function $u$ can be viewed as a difference of two solutions to \eqref{i1} with the same initial datum, hence
they satisfy the bound from (\ref{w1a}), so there exists $m>0$
such that
\begin{equation}\label{w11}
\|u\|_{L_\infty(0,T;L_\infty(\R^n))}\leq m.
\end{equation}
The application of Theorem D to the second integral in  (\ref{w10}) leads us to the following inequality
\begin{equation}\label{w12}
\frac{d}{dt}\|u^2 \|_{L_1(\R^n)} \leq C\| \div\,u\|_{BMO(\R^n)}\|u^2\|_{L_1(\R^n)}\left[
|\ln \|u^2\|_{L_1(\R^n)}| + \ln(e+m) \right]
\end{equation}
with the initial datum $\|u^2|_{t=0}\|_{L_1}=0$ and $m$ from (\ref{w11}).
The Osgood lemma applied to (\ref{w12})  yields the uniqueness to (\ref{i1}).
Note that \eqref{w12} has the form of \eqref{cc4} mentioned in Introduction. Theorem A is proved.

\section{Proof of Theorem B}

The next result concerns the stability of solutions from Theorem A.
We start with the mollified equation \eqref{w3b} for reasons same as previously,
 testing  it  now by $|S_\epsilon(u)|^{p-2}S_\epsilon(u) \pi_r$ with $p$ as in Theorem B.
  Repeating the considerations for (\ref{w4})-(\ref{w10}) we deduce
\begin{equation}\label{b4}
\frac{d}{dt} \int_{\R^n} | u^k-u|^p dx  \leq \int_{\R^n} |\div\,b| \,| u^k-u|^p dx.
\end{equation}

For a given $1\geq \epsilon>0$, we fix $K_\epsilon \in \N$ such that for all $k>K_\epsilon$
\begin{equation}\label{b5}
\|u_0^k-u_0\|_{L_p} \leq \epsilon.
\end{equation}

Let $X_p=|u^k-u|^p$, then by Theorem D ($m$ as in (\ref{w11})) (\ref{b4}) reads
\begin{equation}\label{b6}
\begin{array}{c}
\displaystyle
\frac{d}{dt} \int_{\R^n} X dx \leq \int_{\R^n} |\div\,b| X dx\\[9pt]
\displaystyle
\leq C_0\|\div\,b\|_{BMO(\R^n)}\|X\|_{L_1(\R^n)}\left[ |\ln \|X\|_{L_1(\R^n)}| + \ln (e+2m)\right],
\end{array}
\end{equation}
with $\int_{\R^n} X(x,0) dx \leq \epsilon.$

By  assumptions the r.h.s of (\ref{b6}) is at least  locally integrable, 
hence $\int_{\R^n} X (x,t) dx$ is uniformly continuous. There exits 
a positive time $T_0$
so small that
\begin{equation}\label{b8}
\sup_{t \in [0,T_0]} \int_{\R^n }X(x,t) dx \leq e^{-1}.
\end{equation}
It follows that the function $w|\ln w|$ will be considered as increasing, since $\int_{\R^n} X(x,\cdot)dx$ 
on the chosen time interval takes the values only from the interval $[0,e^{-1}]$. Monotonicity allows us to introduce a function $B: [0,T_0] \to [0,\infty)$ such that
\begin{equation}\label{b9}
\frac{d}{dt} B=C_0 \| \div\,b\|_{BMO(\R^n)} B [|\ln B| + \ln(e+2m)]
\mbox{~~~and~~~} B|_{t=0}=\epsilon.
\end{equation}
The definition of $B$ guarantees that it is an increasing and continuous function, thus there exists $T_1>0$ such that $0< T_1 \leq T_0$ and
\begin{equation}\label{b10}
B(t) \leq e^{-1}< 1 \mbox{~~for~~} t \in [0,T_1].
\end{equation}
Taking the difference between (\ref{b9}) and (\ref{b6}) we get
\begin{equation}\label{b11}
\begin{array}{c}
\displaystyle
 × \frac{d}{dt} [B - \int_{\R^n} X dx]  \geq C_0\| \div b\|_{BMO(\R^n)}\cdot
\\
\displaystyle
\cdot
\left[ B |\ln B| - \int_{\R^n} X dx |\ln \int_{\R^n} X dx| + 
\ln (e+2m) (B -\int_{\R^n} X dx)\right]
\end{array}
\end{equation}
with $B(0)-\int_{\R^n} X(x,0)dx \geq 0$.

Since the monotonicity of the function $w |\ln w|$ on $[0,e^{-1}]$ implies
\begin{equation}\label{b12}
(B |\ln B| - \int_{\R^n} X dx |\ln \int_{\R^n} X dx|)(B -\int_{\R^n} X dx) \geq 0,
\end{equation}
remembering that we consider $t\in [0,T_1]$, from \eqref{b11} we get
\begin{equation}\label{b13}
0\leq \int_{\R^n} X(x,t) dx \leq B(t) \mbox{~~for~~} t\in [0,T_1].
\end{equation}
The above fact reduces our analysis to the considerations of the function $B$. 
Additionally, by the choice of the time interval it follows that  $B(t) <1$ for
$t\in [0,T_1]$, hence
we can   use the estimate
\begin{equation}\label{b14}
|\ln B| \leq \ln \epsilon^{-1}\qquad \mbox{~~for~~} t \in [0,T_1].
\end{equation}
Solving  (\ref{b9}) we get
\begin{equation}\label{b15}
B(t) \leq \epsilon  \, \exp\left\{ C_0[ \ln(e+2m)+\ln \epsilon^{-1}] \int_0^t f(s) ds\right\}
\leq C \epsilon \, \epsilon^{-C_0 \int_0^t f(s)ds},
\end{equation}
where $f(t)=\|\div b(\cdot,t)\|_{BMO(\R^n)}$ and $C$ depends on data given  in Theorems A and B.

Next, we choose $T_2$ so small that $0< T_2 \leq T_1$ and $C_0 \int_0^{T_2} f(s)ds \leq 1/2$, then 
(\ref{b15}) yields
\begin{equation}\label{b16}
\sup_{t \in [0,T_2]} B(t) \leq C \epsilon^{1/2}.
\end{equation}
Here we shall emphasize that $T_2$ is independent from the smallness of $\epsilon$ -- see \eqref{b5}. Thus we are able to start our analysis over
 from the very beginning, but for the initial time $t =T_2$. Since $C_0$ in (\ref{b15}) is an absolute constant we find the next interval $[T_2,T_3]$, 
where  we obtain
\begin{equation}\label{b17}
\sup_{t\in [T_2,T_3]} \|u^k-u\|_{L_p(\R^n)} \leq C\epsilon^{1/4}
\end{equation}
for all $k > K_\epsilon$ -- see (\ref{b5}).
Since $T$ is fixed and finite and by the assumptions $f \in L_1(0,T)$, we are always able to cover the whole interval $[0,T]$ in finite steps, so finally we obtain
\begin{equation}\label{b18}
\sup_{t \in [0,T]} B(t) \leq C \epsilon^{a}.
\end{equation}
with $a>0$ defined by the properties of $f$ and again $C$ depending on all data, but  independent from $\epsilon$. 
Letting $\epsilon \to 0$ we prove (\ref{b3}). Theorem B is proved.

\section{Proof of Theorem C}

Our last result describes the uniqueness criteria for weak solutions, provided their
existence in the $L_\infty(0,T;L_p(\R^n))$-class.
 The problem reduces to (\ref{i1}) with zero initial data and $u \in L_\infty(0,T;L_p(\R^n))$. 
To work in optimal regularity of coefficients we consider (\ref{w4})
\begin{equation*}\label{c1}
 \d_t S_\epsilon(u)+b \cdot \nabla S_\epsilon(u)= R_\epsilon \to 0 
\mbox{~~in~~} L_1(0,T;L_{1(loc)}(\R^n)).
\end{equation*}
Next, we introduce the renormalized solution for (\ref{i1}) -- we refere here to \cite{DL} where this approach has been developed. Take $\beta \in C^1(\R)$, 
i.e. $\|\beta\|_{L_\infty(\R)}+\|\beta'\|_{L_\infty(\R)} < \infty$, then
\begin{equation}\label{c2}
 \d_t \beta(S_\epsilon(u))+b \cdot \nabla \beta(S_\epsilon(u))= R_\epsilon \beta'(S_\epsilon(u))
\end{equation}
which implies the limit for $\epsilon \to 0$
\begin{equation}\label{c3}
 \d_t \beta(u)+b\cdot \nabla \beta(u)=0.
\end{equation}
As the function $\beta$ we choose $T_m: \R \to [0,m^p]$ such that
\begin{equation}\label{c4}
 T_m(s)=\left\{ 
\begin{array}{lcr}
 × |s|^p & \mbox{for} & |s|< m \\
m^p & \mbox{for} & |s| \geq m
\end{array}
\right.
\end{equation}
defined for fixed $m \in \R_+$. $T_m$ is not a $C^1$-function, but a simple approximation procedure will lead us to \eqref{c3} with $\beta=T_m$.

Since we do not control  integrability of all terms in (\ref{c3}) we use the function $\pi_r$ from 
(\ref{w6}) to localize the problem getting
\begin{equation}\label{c5}
 \frac{d}{dt} \int_{\R^n} T_m(u)\pi_r dx \leq \int_{\R^n} | \div b | \, T_m(u) \pi_r dx + 
\int_{\R^n} | b \cdot \nabla \pi_r| \, T_m(u) dx.
\end{equation}
For fixed $m$ and $r$ letting to infinity the last term vanishes, so we obtain
\begin{equation}\label{c6}
 \frac{d}{dt} \int_{\R^n} T_m(u) dx \leq \int_{\R^n} | \div b | \, T_m(u)  dx.
\end{equation}
Theorem D applied to the r.h.s. of (\ref{c6}) yields
\begin{equation}\label{c7}
 \frac{d}{dt} \| T_m(u)\|_{L_1(\R^n)} \leq C\| \div b \|_{BMO(\R^n)} 
\|T_m(u)\|_{L_1(\R^n)}\left[ |\ln \|T_m(u)\|_{L_1(\R^n)}|
+\ln (e+m^p)\right]
\end{equation}
with $\|T_m(u)(\cdot,0)\|_{L_1(\R^n)}=0$.

The same way as in the proof of Theorem A, the Osgood lemma yields $T_m(u)\equiv 0$. Letting $m \to \infty$, by (\ref{c4}) we conclude $u\equiv 0$.
Thus, $u_1\equiv u_2$. Theorem C is proved.

\section{Appendix}

{\it Sketch of the proof of Theorem D.} We proceed almost as in \cite{MR}. The assumption of the boundedness of supp $f$ allows us to consider the studied integral (the l.h.s. of \eqref{d1}) on a torus $\T^n=\R^n/ (d\Z^n) (=[0,d)^n)$ with sufficiently large $d$ guaranteeing that supp $f$ can be treated as a subset of $\T^n$.
Consider  the Hardy space on $\T^n$ with the  following norm
\begin{equation}\label{x0}
\|g\|_{{\cal H}^1(\T^n)}=\|g\|_{L_1(\T^n)}+\sum_{k=1}^n \|R_kg\|_{L_1(\T^n)},
\end{equation}
where $R_k$ are the Riesz operators -- \cite{S},\cite{Tor}. Since $BMO(\T^n)=({\cal H}^1(\T^n))^*$, we get
\begin{equation}\label{x1}
|\int_{\T^n} fg dx | \leq \|f\|_{BMO(\T^n)}\|g\|_{{\cal H}^1(\T^n)}.
\end{equation}
Hence to control the norm  (\ref{x0}) 
 an estimate of $\|R_kg\|_{L_1(\T^n)}$ is required. The classical Zygmund's result \cite{Z} says:
\begin{equation}\label{x1a}
\|R_k h \|_{L_1(\T^n)} \leq C+C \int_{\T^n}  |h| \ln^+|h| dx,
\end{equation}
where $\ln^+ a = \max\{ \ln a ,0\}$ and constants $C$ depends on $d$, so on the diameter of supp $f$.

Let us observe that $\ln^+(g/\lambda)=\ln g - \ln \lambda$ for $g\geq \lambda$ and
\begin{equation}\label{x4}
 |\ln g |_{g\geq \lambda}|\leq \ln(1+\|g\|_{L_\infty(\T^n)})+|\ln \frac{g}{1+\|g\|_{L_\infty(\T^n)}}|_{g\geq \lambda}|
\leq 2 \ln(1+\|g\|_{L_\infty(\T^n)})+|\ln \lambda |.
\end{equation}
Taking $h=\frac{g}{\|g\|_{L_1(\T^n)}}$ in \eqref{x1a}, employing  \eqref{x4}, we conclude
\begin{equation}\label{x5}
 \| R_kg\|_{L_1(\T^n)} \leq C\|g\|_{L_1(\T^n)}+C\int_{\T^n} |g| \left[\ln (1+\|g\|_{L_\infty(\T^n)})+|\ln \|g\|_{L_1(\T^n)}|
\right]dx
\end{equation}
Inequalities \eqref{x1}, \eqref{x5} yields \eqref{d1}.

\smallskip

{\it The  commutator estimate.} Let us recall the well known facts concerning the mollification of the equation and the behavior of the commutators. 

Introduce $m_1: \R^n \to [0,\infty)$ such that
$$
m_1(x)= N_n \left\{
\begin{array}{lcr}
\exp\{ - \frac{1}{1-|x|^2} \} & \mbox{for} & |x|  < 1\\
0 & \mbox{for} & |x| \geq 1
\end{array}
\right.
$$
where  the number $N_n$ is determined by the constraint $\int_{\R^n} m_1 dx =1$. Then 
for given $\epsilon >0$
we define
\begin{equation}\label{x1b}
m_\epsilon(x):= \frac{1}{\epsilon^n} m_1(\frac x\epsilon) \qquad 
\mbox{with }  \int_{\R^n} m_\epsilon dx =1.
\end{equation}
It is clear that $m_\epsilon \to \delta$ in ${\cal D}'(\R^n)$,
where $\delta$ is the Dirac mass located at the origin of $\R^n$. The function $m_\epsilon$ introduces  an operator $S_\epsilon: L_{1(loc)}(\R^n) \to C^\infty(\R^n)$ 
\begin{equation}\label{x1c}
 S_\epsilon(h)=m_\epsilon \ast h=\int_{\R^n} m_\epsilon(x-y)h(y) dy.
\end{equation}
The standard theory,  see \cite{DL}, guarantees the following estimate for the 
commutator
\begin{equation}\label{x2}
 b\cdot \nabla S_\epsilon(u)- S_\epsilon(b \cdot \nabla u) \to 0 
\mbox{~~in~~} L_1(0,T;L_{1(loc)}(\R^n)),
\end{equation}
provided that $u \in L_\infty(0,T;L_{p(loc)}(\R^n))$, 
$b \in L_1(0,T;W^1_{p'(loc)}(\R^n))$ and $1=\frac{1}{p}+\frac{1}{p'}$ for
$p\in [1,\infty]$.
The proof of (\ref{x2}) belongs to the by now  classical theory, we omit it here and refere 
again to \cite{DL}.

\medskip

\noindent
{\footnotesize {\bf Acknowledgments.} The author thanks  Beno{\^{\i}}t Perthame and Walter Rusin for fruitful discussion.
The  author has been supported by
Polish  grant No. 1 P03A 021 30 and by ECFP6 M.Curie ToK program SPADE2,
MTKD-CT-2004-014508 and SPB-M. He also thanks Laboratoire Jacques-Louis Lions and Leibniz Universit\"at Hannover, where 
parts of this paper were performed, for their hospitality. The stay in Hannover has been supported by the Humboldt Foundation.}

\end{document}